\documentclass{tran-l}
\usepackage{graphicx}

 \newtheorem{thm}{Theorem}[section]
 \newtheorem{cor}[thm]{Corollary}
 \newtheorem{ass}[thm]{Assumption}
 \newtheorem{lem}[thm]{Lemma}
 
 \theoremstyle{definition}
 
 \theoremstyle{remark}
 \newtheorem{rem}[thm]{Remark}
 \numberwithin{equation}{section}

\begin{document}

\title[Block-Diagonalization of the Linearized Coupled-Mode System]{Block-Diagonalization of the Linearized Coupled-Mode System}

\author{Marina Chugunova,  Dmitry Pelinovsky}

\address{Department of Mathematics, McMaster University, Hamilton, Ontario, L8S 4K1, Canada}

\email{chugunom@math.mcmaster.ca; dmpeli@math.mcmaster.ca}

\thanks{This work was completed with the support of the SharcNet Graduate Scholarship.}


\keywords{Hamiltonian coupled-mode systems, gap solitons, spectral
stability, invariant subspaces, eigenvalues}


\begin{abstract}
We consider the Hamiltonian coupled-mode system that occur in
nonlinear optics, photonics, and atomic physics. Spectral stability
of gap solitons is determined by eigenvalues of the linearized
coupled-mode system, which is equivalent to a four-by-four Dirac
system with sign-indefinite metric. In the special class of
symmetric nonlinear potentials, we construct a block-diagonal
representation of the linearized equations, when the spectral
problem reduces to two coupled two-by-two Dirac systems. The
block-diagonalization is used in numerical computations of
eigenvalues that determine stability of gap solitons.
\end{abstract}

\maketitle

\section{Introduction}

Various applications in nonlinear optics \cite{SS94}, photonics
band-gap engineering \cite{photonics} and atomic physics
\cite{atomic} call for systematic studies of the {\em coupled-mode
system}, which is expressed by two first-order semi-linear PDEs in
one space and one time dimensions. In nonlinear optics, the
coupled-mode system describes counter-propagating light waves, which
interact with a linear grating in an optical waveguide \cite{S98}.
In photonics, the coupled-mode system is derived for coupled
resonant waves in stop bands of a low-contrast three-dimensional
photonic crystal \cite{AP05}. In atomic physics, the coupled-mode
system describes matter-wave Bose-Einstein condensates trapped in an
optical lattice \cite{PSK04}. Existence, stability and nonlinear
dynamics of {\em gap solitons}, which are localized solutions of the
coupled-mode system, are fundamental problems for interest in the
aforementioned physical disciplines.

In the context of spectral stability of gap solitons, it has been
discovered that the linearized coupled-mode system is equivalent to
a four-by-four Dirac system with sign-indefinite metric, where
numerical computations of eigenvalues represent a difficult
numerical task. The pioneer work in \cite{B98,B99} showed that
spurious unstable eigenvalues originate from the continuous spectrum
in the Fourier basis decomposition and the Galerkin approximation. A
delicate but time-consuming implementation of the continuous Newton
method was developed to identify true unstable eigenvalues from the
spurious ones \cite{B98,B99}. Similar problems were discovered in
the variational method \cite{KL96,KL97} and in the numerical
finite-difference method \cite{M99,M00}.

While some conclusions on instability bifurcations of gap solitons
in the coupled-mode equations can be drawn on the basis of
perturbation theory \cite{B98} and Evans function methods
\cite{KS02,PS03}, the numerical approximation of eigenvalues was an
open problem until recently. A new progress was made with the use of
exterior algebra in the numerical computations of the Evans function
\cite{DG04}, when the same results on instability bifurcations of
gap solitons as in \cite{B98} were recovered. Similar shooting
method was also applied to gap solitons in a more general model of a
nonlinear Schr\"{o}dinger equation with a periodic potential
\cite{PSK04}.

Our work addresses the problem of numerical approximations of
eigenvalues of the linearized coupled-mode system with a different
objective. We will show that the linearized coupled-mode system with
a symmetric potential function can be block-diagonalized into two
coupled two-by-two Dirac systems. The two Dirac systems represent
the linearized Hamiltonian of the coupled-mode equations and
determine instability bifurcations and unstable eigenvalues of gap
solitons.

The purpose of block-diagonalization is twofold. First, the number
of unstable eigenvalues and details of instability bifurcations can
be investigated analytically from the number of non-zero isolated
eigenvalues of the linearized Hamiltonian. This analysis will be
reported in the forthcoming publication. Second, a numerical
algorithm can be developed to compute efficiently the entire
spectrum of the linearized coupled-mode system. These numerical
results are reported here for an example of symmetric quadric
potential functions.

The paper is organized as follows. Section 2 describes the model and
its symmetries. Section 3 gives construction and properties of gap
solitons in the nonlinear coupled-mode system. Section 4 presents
block-diagonalization of the linearized coupled-mode system. Section
5 contains numerical computations of the spectrum of the
block-diagonalized system. Appendix A presents derivation of exact
solutions for gap solitons in the coupled-mode system with symmetric
quadric potential functions.

\section{Coupled-mode system}

We consider the Hamiltonian coupled-mode system in the form:
\begin{eqnarray}
\label{coupled-mode-system} \left \{
 \begin{array}{c}
 i(u_t + u_x) + v = \partial_{\bar{u}} W(u,\bar{u},v,\bar{v}) \\
 i(v_t - v_x) + u = \partial_{\bar{v}} W(u,\bar{u},v,\bar{v})
 \end{array}
 \right.
\end{eqnarray}
where $(u,v) \in \mathbb{C}^2$, $x \in \mathbb{R}$, $t \geq 0$, and
$W(u,\bar{u},v,\bar{v})$ is real-valued. We assume that the
potential function satisfies the following three conditions:

\begin{enumerate}
\item $W$ is invariant with respect to the gauge transformation:
$(u,v) \mapsto e^{i \alpha} (u,v)$, for all $\alpha \in \mathbb{R}$

\item $W$ is symmetric with respect to the interchange:
$(u,v) \mapsto (v,u)$

\item $W$ is analytic in its variables near $u = v = 0$,
such that $W = O(4)$.
\end{enumerate}

The first property is justified by the standard derivation of the
coupled-mode system (\ref{coupled-mode-system}) with an envelope
approximation \cite{AP05}. The second property defines a class of
symmetric nonlinear potentials. Although it is somewhat restrictive,
symmetric nonlinear potentials are commonly met in physical
applications of the system (\ref{coupled-mode-system}). The third
property is related to the normal form analysis \cite{SU03}, where
the nonlinear functions are approximated by Taylor polynomials.
Since the quadratic part of the potential function is written in the
left-hand-side of the system (\ref{coupled-mode-system}) and the
cubic part violates the gauge transformation and analyticity
assumptions, the Taylor polynomials of $W$ start with quadric terms,
denoted as $O(4)$.

We find a general representation of the function
$W(u,\bar{u},v,\bar{v})$ that satisfies the conditions (1)-(3) and
list all possible (four-parameter) quadric terms of $W$.

\begin{lem}
If $W \in \mathbb{C}$ and property (1) is satisfied, such that
\begin{equation}
\label{gauge-identity} W(u,\bar{u},v,\bar{v}) = W\left(u e^{i
\alpha}, \bar{u} e^{-i \alpha}, v e^{i \alpha}, \bar{v} e^{- i
\alpha} \right), \qquad \forall \alpha \in \mathbb{R},
\end{equation}
then $W = W(|u|^2,|v|^2,u\bar{v})$.
\end{lem}
\begin{proof}
By differentiating (\ref{gauge-identity}) in $\alpha$ and setting
$\alpha = 0$, we have the differential identity:
\begin{equation}
\label{PDO} D W \equiv i \left( u\frac{\partial}{\partial u} -
\bar{u}\frac{\partial}{\partial \bar{u}} + v
\frac{\partial}{\partial v} - \bar{v}\frac{\partial}{\partial
\bar{v}} \right) W(u,\bar{u},v,\bar{v}) = 0.
\end{equation}
Consider the set of quadratic variables
$$
z_1 = |u|^2, \quad z_2 = |v|^2, \quad z_3 = \bar{u} v, \quad z_4 =
u^2,
$$
which is independent for any $u \neq 0$ and $v \neq 0$ in the sense
that the Jacobian is non-zero. It is clear that $D z_{1,2,3} = 0$
and $D z_4 = 2 z_4$.  Therefore, $D W = 2 z_4 \partial_{z_4} W = 0$,
such that $W = W(z_1,z_2,z_3)$.
\end{proof}

\begin{cor}
\label{cor22} If $W \in \mathbb{R}$ and property (1) is met, then $W
= W(|u|^2,|v|^2,u\bar{v}+v\bar{u})$.
\end{cor}

\begin{lem}
If $W \in \mathbb{R}$ and properies (1)-(3) are satisfied, then $W =
W(|u|^2 + |v|^2,|u|^2 |v|^2, u \bar{v} + v \bar{u})$.
\end{lem}

\begin{proof}
By Corollary \ref{cor22} and property (2), we can re-order the
arguments of $W$ as $W = W(|u| + |v|,|u||v|,u \bar{v} + v \bar{u})$.
By analyticity in property (3), $W$ may depend only on $|u|^2$ and
$|v|^2$ rather than on $|u|$ and $|v|$.
\end{proof}

\begin{cor}
\label{Re0} If $W \in \mathbb{R}$ and properties (1)-(3) are
satisfied, then
\begin{equation}
\label{W-property1} \left( u \frac{\partial}{\partial u} +
\bar{u}\frac{\partial}{\partial \bar{u}} - v
\frac{\partial}{\partial v} - \bar{v}\frac{\partial}{\partial
\bar{v}} \right) W(u,\bar{u},v,\bar{v}) \biggr|_{|u|^2 = |v|^2} = 0
\end{equation}
\end{cor}

\begin{cor}
The only quadric potential function $W \in \mathbb{R}$ that
satisfies properties (1)-(3) is given by
\begin{equation}
\label{example} W = \frac{a_1}{2} (|u|^4 + |v|^4) + a_2 |u|^2 |v|^2
+ a_3 (|u|^2+|v|^2)(v \bar{u} + \bar{v} u) + \frac{a_4}{2}(v \bar{u}
+ \bar{v} u)^2,
\end{equation}
where $(a_1,a_2,a_3,a_4)$ are real-valued parameters. It follows
then that
\begin{eqnarray}
 \left \{
 \begin{array}{c}
 \partial_{\overline{u}} W = a_1 |u|^2 u + a_2 u |v|^2 + a_3
 \left[ (2 |u|^2 + |v|^2)v + u^2 \bar{v} \right] + a_4 \left[
 v^2 \bar{u} + |v|^2 u \right] \nonumber \\
 \partial_{\overline{v}} W = a_1 |v|^2 v + a_2 v |u|^2 + a_3
 \left[ (2 |v|^2 + |u|^2) u + v^2 \bar{u} \right] + a_4 \left[
 u^2 \bar{v} + |u|^2 v \right] \nonumber
 \end{array}
 \right.
\end{eqnarray}
\end{cor}

The potential function (\ref{example}) with $a_1,a_2 \neq 0$ and
$a_3 = a_4 = 0$ represents a standard coupled-mode system for a
sub-harmonic resonance, e.g. in the context of optical gratings with
constant Kerr nonlinearity \cite{SS94}. When $a_1 = a_3 = a_4 = 0$,
this system is integrable with inverse scattering and is referred to
as the massive Thirring model \cite{MTM}. When $a_1 = a_2 = 0$ and
$a_3,a_4 \neq 0$, the coupled-mode system corresponds to an optical
grating with varying, mean-zero Kerr nonlinearity, where $a_3$ is
the Fourier coefficient of the resonant sub-harmonic and $a_4$ is
the Fourier coefficient of the non-resonant harmonic \cite{AP05}
(see also \cite{S98}).

We rewrite the coupled-mode system (\ref{coupled-mode-system}) as a
Hamiltonian system in complex-valued matrix-vector notations:
\begin{equation}
\label{Ham-system} \frac{d {\bf u}}{d t} = J \nabla H({\bf u}),
\end{equation}
where ${\bf u} = (u,\bar{u},v,\bar{v})^T$,
$$
J = \left [ \begin{array}{cccc} 0 & -i & 0 & 0 \\
i & 0 & 0 & 0 \\ 0 & 0 & 0 & -i \\ 0 & 0 & i & 0  \end{array}
\right] = -J^T,
$$
and $H(u,\bar{u},v,\bar{v}) = \int_{\mathbb{R}}
h(u,\bar{u},v,\bar{v}) dx$ is the Hamiltonian functional with the
density:
$$
h = W(u, \bar{u}, v, \bar{v}) - (v \bar{u} + u \bar{v}) +
\frac{i}{2} (u \bar{u}_x - u_x \bar{u}) - \frac{i}{2}(v \bar{v}_x -
v_x \bar{v}).
$$
The Hamiltonian $H(u,\bar{u},v,\bar{v})$ is constant in time $t \geq
0$. Due to the gauge invariance, the coupled-mode system
(\ref{coupled-mode-system}) has another constant of motion
$Q(u,\bar{u},v,\bar{v})$, where
\begin{equation}
\label{quantity-Q} Q = \int_{\mathbb{R}} \left( |u|^2+|v|^2 \right)
dx.
\end{equation}
Conservation of $Q$ can be checked by direct computation:
\begin{equation}
\label{balance-equation} \frac{\partial}{\partial t}(|u|^2 + |v|^2)
+ \frac{\partial}{\partial x}(|u|^2 - |v|^2) = D W = 0,
\end{equation}
where the operator $D$ is defined in (\ref{PDO}). Due to the
translational invariance, the coupled-mode system
(\ref{coupled-mode-system}) has yet another constant of motion
$P(u,\bar{u},v,\bar{v})$, where
\begin{equation}
\label{quantity-P} P = \frac{i}{2} \int_{\mathbb{R}} \left( u
\bar{u}_x - u_x \bar{u} + v \bar{v}_x - v_x \bar{v} \right) dx.
\end{equation}
In applications, the quantities $Q$ and $P$ are referred to as the
power and momentum of the coupled-mode system.

\section{Existence of gap solitons}

{\em Stationary} solutions of the coupled-mode system
(\ref{coupled-mode-system}) take the form:
\begin{equation}
\label{soliton} \left \{
\begin{array}{c}
u_{\rm st}(x,t) = u_0(x+s) e^{i\omega t + i \theta}\\
v_{\rm st}(x,t) = v_0(x+s) e^{i\omega t + i \theta}
\end{array}
\right.
\end{equation}
where $(s,\theta) \in \mathbb{R}^2$ are arbitrary parameters, while
the solution $(u_0,v_0) \in \mathbb{C}^2$ on $x \in \mathbb{R}$ and
the domain for parameter $\omega \in \mathbb{R}$ are to be found
from the nonlinear ODE system:
\begin{equation}
\label{ODE-system} \left \{
\begin{array}{c}
iu_0' = \omega u_0 - v_0 + \partial_{\bar{u}_0}
W(u_0,\bar{u}_0,v_0,\bar{v}_0) \\
-iv_0' = \omega v_0 - u_0 + \partial_{\bar{v}_0}
W(u_0,\bar{u}_0,v_0,\bar{v}_0)
\end{array} \right.
\end{equation}
Stationary solutions are critical points of the Lyapunov functional:
\begin{equation}
\label{Lyapunov-function} \Lambda = H(u,\bar{u},v,\bar{v}) + \omega
Q(u,\bar{u},v,\bar{v}),
\end{equation}
such that variations of $\Lambda$ produce the nonlinear ODE system
(\ref{ODE-system}).

\begin{lem}
\label{lem-u0-v0} Assume that there exists a decaying solution
$(u_0,v_0)$ of the system (\ref{ODE-system}) on $x \in \mathbb{R}$.
If $W \in \mathbb{R}$ satisfies properties (1)-(3), then $u_0 =
\bar{v}_0$ (module to an arbitrary phase).
\end{lem}
\begin{proof}
It follows from the balance equation (\ref{balance-equation}) for
the stationary solutions (\ref{soliton}) that
$$
|u_0|^2 - |v_0|^2 = C_0 = 0, \qquad \forall x \in \mathbb{R},
$$
where the constant $C_0 = 0$ is found from decaying conditions at
infinity. Let us represent the solutions $(u_0,v_0)$ in the form:
\begin{equation}
\label{representation-soliton}
 \left \{
\begin{array}{c}
u_0(x)= \sqrt{Q(x)} e^{i \Theta(x) + i \Phi(x)}\\
v_0(x)=\sqrt{Q(x)} e^{- i \Theta(x) + i \Phi(x)}
\end{array}
 \right.
\end{equation}
such that
\begin{equation}
\label{ODE-system-1} \left \{
\begin{array}{c}
i Q' - 2Q (\Theta' + \Phi') = 2\omega Q - 2 Q e^{-2i\Theta} + 2
\bar{u}_0 \partial_{\bar{u}_0} W(u_0,\bar{u}_0,v_0,\bar{v}_0) \\
-i Q' - 2Q (\Theta' - \Phi') = 2\omega Q - 2 Q e^{2i\Theta} + 2
 \bar{v}_0 \partial_{\bar{v}_0} W(u_0,\bar{u}_0,v_0,\bar{v}_0)
\end{array}
\right.
\end{equation}
Separating the real parts, we obtain
\begin{equation}
\left \{
\begin{array}{c}
Q ( \cos(2 \Theta) - \omega - \Theta' - \Phi') = {\rm Re}
\left[ \bar{u}_0 \partial_{\bar{u}_0} W(u_0,\bar{u}_0,v_0,\bar{v}_0) \right] \\
Q ( \cos(2 \Theta) - \omega - \Theta' + \Phi') = {\rm Re} \left[
\bar{v}_0 \partial_{\bar{v}_0} W(u_0,\bar{u}_0,v_0,\bar{v}_0)
\right]
\end{array}
\right.
\end{equation}
By Corollary \ref{Re0}, we have $\Phi' \equiv 0$, such that $\Phi(x)
= \Phi_0$.
\end{proof}

\begin{cor}
Let $u_0 = \bar{v}_0$. The ODE system (\ref{ODE-system}) reduces to
the planar Hamiltonian form:
\begin{equation}
\label{ODE-Hamiltonian} \frac{d}{dx} \left( \begin{array}{cc} p \\
q \end{array} \right) = \left( \begin{array}{ccc} 0 & -1 \\ +1 & 0
\end{array} \right) \nabla h(p,q),
\end{equation}
where $p = 2 \Theta$, $q = Q$, and
\begin{equation}
\label{Hamiltonian-function-ODE} h = \tilde{W}(p,q) - 2 q \cos p + 2
\omega q, \qquad \tilde{W}(p,q) = W(u_0,\bar{u}_0,v_0,\bar{v}_0).
\end{equation}
\end{cor}

\begin{proof}
In variables $(Q,\Theta)$ defined by (\ref{representation-soliton})
with $\Phi(x) = \Phi_0 = 0$, we rewrite the ODE system
(\ref{ODE-system-1}) as follows:
\begin{equation}
\label{ODE-system-2} \left \{
\begin{array}{c}
Q' = 2 Q \sin(2 \Theta) + 2 {\rm Im} \left[
\bar{u}_0 \partial_{\bar{u}_0} W(u_0,\bar{u}_0,v_0,\bar{v}_0) \right] \\
Q \Theta' = -\omega Q + Q \cos(2 \Theta) - {\rm Re} \left[
 \bar{u}_0 \partial_{\bar{u}_0} W(u_0,\bar{u}_0,v_0,\bar{v}_0)
 \right]
\end{array}
\right.
\end{equation}
The system (\ref{ODE-system-2}) is equivalent to the Hamiltonian
system (\ref{ODE-Hamiltonian}) and (\ref{Hamiltonian-function-ODE})
if
\begin{equation}
\label{ODE-system-2-too} \left \{
\begin{array}{c}
\partial_p \tilde{W}(p,q) = i \left[ u_0 \partial_{u_0}
- \bar{u}_0 \partial_{\bar{u}_0} \right] W(u_0,\bar{u}_0,v_0,\bar{v}_0) \\
q \partial_q \tilde{W}(p,q) =  \left[  u_0 \partial_{u_0} +
\bar{u}_0 \partial_{\bar{u}_0} \right]
W(u_0,\bar{u}_0,v_0,\bar{v}_0)
\end{array}
\right.
\end{equation}
The latter equations follows from (\ref{PDO}), (\ref{W-property1}),
and (\ref{representation-soliton}) with the chain rule.
\end{proof}

\begin{cor}
\label{cor-property2} Let $u_0 = \bar{v}_0$. Then,
\begin{equation}
\label{W-property2}
\partial^2_{u_0 \bar{u}_0} W = \partial^2_{v_0 \bar{v}_0} W, \quad
\partial^2_{\bar{u}_0^2} W = \partial^2_{v_0^2} W, \quad
\partial^2_{u_0 v_0} W = \partial^2_{\bar{u}_0 \bar{v}_0} W.
\end{equation}
\end{cor}

\begin{rem}
The family of stationary solutions (\ref{soliton}) can be extended
to the family of travelling solutions of the coupled-mode system
(\ref{coupled-mode-system}) by means of the Lorentz transformation
\cite{DG04}. With the boosted variables,
$$
X = \frac{x-ct}{\sqrt{1 - c^2}}, \quad T = \frac{t -
cx}{\sqrt{1-c^2}}, \quad U = \left(\frac{1-c}{1+c}\right)^{1/4} u,
\quad V = \left(\frac{1+c}{1-c}\right)^{1/4} v,
$$
where $c \in (-1,1)$, the family of travelling solutions still
satisfies the constraint $|U_0|^2 = |V_0|^2$ from the balance
equation (\ref{balance-equation}). However, Corollary \ref{Re0}
fails for a boosted potential function
$\tilde{W}(U,\bar{U},V,\bar{V})$ and the representation
(\ref{representation-soliton}) results no longer in the relation
$U_0 = \bar{V}_0$ \cite{B98}. It will be studied separately if the
block-diagonalization of the linearized coupled-mode system can be
extended (in a non-trivial matter) to the family of travelling
solutions.
\end{rem}

Decaying solutions of the system (\ref{ODE-system}) with a
homogeneous polynomial function $W(u,\bar{u},v,\bar{v})$ are
analyzed in Appendix A. Conditions for their existence are
identified for the quadratic potential function (\ref{example}).
Decaying solutions may exist in the gap of continuous spectrum of
the coupled-mode system (\ref{coupled-mode-system}) for $\omega \in
(-1,1)$. We introduce two auxiliary parameters:
\begin{equation}
\label{parameter-ODEs} \mu = \frac{1 - \omega}{1 + \omega}, \qquad
\beta = \sqrt{1 - \omega ^2},
\end{equation}
such that $0 < \mu < \infty$ and $0 < \beta \leq 1$. When $a_1 = 1$,
$a_2 = \rho$, and $a_3 = a_4 = 0$, we obtain in Appendix A the
decaying solution $u_0(x)$ in the explicit form:
\begin{equation}
\label{soliton-explicit-1} u_0 =
\sqrt{\frac{2(1-\omega)}{1+\rho}}\frac{1}{(\cosh{\beta x} + i
\sqrt{\mu} \sinh{\beta x})}.
\end{equation}
When $\omega \to 1$ (such that $\mu \to 0$ and $\beta \to 0$), the
decaying solution (\ref{soliton-explicit-1}) becomes small in
absolute value and approaches the limit of ${\rm sech}$-solutions
${\rm sech}(\beta x)$. When $\omega \to -1$ (such that $\mu \to
\infty$ and $\beta \to 0$), the decaying solution
(\ref{soliton-explicit-1}) remains finite in absolute value and
approaches the limit of the algebraically decaying solution:
$$
u_0 = \frac{2}{\sqrt{1 + \rho} (1 + 2 i x)}.
$$
When $a_1 = a_2 =0, a_3 = 1$ and $a_4 = s$, the decaying solution
$u_0(x)$ exists in two sub-domains: $\omega > 0$, $s > -1$ and
$\omega < 0$, $s < 1$. When $\omega > 0$, $s > -1$, the solution
takes the form:
\begin{equation}
\label{soliton-explicit-2a} u_0 = \sqrt{\frac{1 - \omega}{2}}
\frac{(\cosh{\beta x} - i \sqrt {\mu} \sinh {\beta
x})}{\sqrt{\Delta_+(x)}},
\end{equation}
where
$$
\Delta_+ = [(s-1) \mu^2- 2 s \mu + (s + 1)]\cosh^4(\beta x) + 2 [s
\mu - (s - 1) \mu^2] \cosh^2(\beta x) + (s-1)\mu^2.
$$
When $\omega < 0$, $s < 1$, the solution takes the form:
\begin{equation}
\label{soliton-explicit-2b} u_0 = \sqrt{\frac{1 - \omega}{2}}
\frac{(\sinh {\beta x} - i \sqrt {\mu} \cosh {\beta
x})}{\sqrt{\Delta_-(x)}}.
\end{equation}
where
$$
\Delta_- = [(s+1) -2 s \mu - (s-1) \mu^2]\cosh^4(\beta x) + 2 [s + 1
- s \mu] \cosh^2(\beta x) - (s+1).
$$
In both limits $\omega \to 1$ and $\omega \to -1$, the decaying
solutions (\ref{soliton-explicit-2a}) and
(\ref{soliton-explicit-2b}) approach the small-amplitude ${\rm
sech}$-solution ${\rm sech}(\beta x)$. In the limit $\omega \to 0$,
the decaying solutions (\ref{soliton-explicit-2a}) and
(\ref{soliton-explicit-2b}) degenerate into a non-decaying bounded
solution with $|u_0(x)|^2 = \frac{1}{2}$.

\section{Block-diagonalization of the linearized system}

Linearization of the coupled-mode system (\ref{coupled-mode-system})
at the stationary solutions (\ref{soliton}) with $s = \theta = 0$ is
defined as follows:
\begin{equation}
\label{linearization} \left \{
\begin{array}{c}
u(x,t) = e^{i\omega t} \left( u_0(x) + U_1(x) e^{\lambda t} \right) \\
\bar{u}(x,t) = e^{-i\omega t} \left( \bar{u}_0(x) + U_2(x) e^{\lambda t} \right) \\
v(x,t) = e^{i\omega t} \left( v_0(x)+ U_3(x) e^{\lambda t} \right) \\
\bar{v}(x,t) = e^{-i\omega t} \left( \bar{v}_0(x)+ U_4(x)
e^{\lambda t} \right)
\end{array}
\right.
\end{equation}
where $v_0 = \bar{u}_0$, according to Lemma \ref{lem-u0-v0}. Let
$({\bf f},{\bf g})$ be a standard inner product for ${\bf f},{\bf g}
\in L^2(\mathbb{R},\mathbb{C}^4)$. Expanding the Lyapunov functional
(\ref{Lyapunov-function}) into Taylor series near ${\bf u}_0 =
(u_0,\bar{u}_0,v_0,\bar{v}_0)^T$, we have:
\begin{equation}
\Lambda = \Lambda({\bf u_0}) + \left( {\bf U}, \nabla \Lambda|_{{\bf
u}_0} \right) + \frac{1}{2} \left( {\bf U}, H_{\omega} {\bf U}
\right) + \ldots,
\end{equation}
where ${\bf U} = (U_1,U_2,U_3,U_4)^T$ and $H_{\omega}$ is the the
linearized energy operator in the explicit form
\begin{equation}
H_{\omega} = D(\partial_x) + V(x),
\end{equation}
where
\begin{equation}
\label{operator-D} D = \left( \begin{array}{cccc} \omega - i
\partial_x & 0 & -1 & 0 \\ 0 & \omega + i \partial_x & 0 & -1 \\
-1 & 0 & \omega + i \partial_x & 0 \\
0 & -1 & 0 & \omega - i \partial_x \end{array}\right)
\end{equation}
and
\begin{equation}
\label{operator-V} V = \left(
\begin{array}{cccc} \partial^2_{\bar{u}_0 u_0} &
\partial^2_{\bar{u}_0^2} &  \partial^2_{\bar{u}_0 v_0} &  \partial^2_{\bar{u}_0 \bar{v}_0} \\
\partial^2_{u_0^2} &  \partial^2_{u_0 \bar{u}_0} &  \partial^2_{u_0 v_0} &  \partial^2_{u_0 \bar{v}_0} \\
\partial^2_{\bar{v}_0 u_0} &  \partial^2_{\bar{v}_0 \bar{u}_0} &  \partial^2_{\bar{v}_0 v_0} &  \partial^2_{\bar{v}_0^2}\\
\partial^2_{v_0 u_0} &  \partial^2_{v_0 \bar{u}_0} &  \partial^2_{v_0^2} &  \partial^2_{v_0 \bar{v}_0}
\end{array}\right) W(u_0,\bar{u}_0,v_0,\bar{v}_0).
\end{equation}
The linearization (\ref{linearization}) of the nonlinear
coupled-mode system (\ref{coupled-mode-system}) results in the
linearized coupled-mode system in the form:
\begin{equation}
\label{linearized-Ham-system} H_{\omega} {\bf U} = i \lambda \sigma
{\bf U},
\end{equation}
where $\sigma$ is a diagonal matrix of $(1,-1,1,-1)$. Due to the
gauge and translational symmetries, the energy operator
$H_{\omega}$ has a non-empty kernel which includes two
eigenvectors:
\begin{equation}
\label{kernel} {\bf U}_1 = \sigma {\bf u}_0(x), \qquad {\bf U}_2 =
{\bf u}_0'(x).
\end{equation}
The eigenvectors ${\bf U}_{1,2}$ represent derivatives of the
stationary solutions (\ref{soliton}) with respect to parameters
$(\theta,s)$. We adopt a standard assumption that the coupled-mode
system is generic.

\begin{ass}
The kernel of $H_{\omega}$ is exactly two-dimensional with the
eigenvectors (\ref{kernel}).
\end{ass}

Due to the Hamiltonian structure, the linearized operator $\sigma
H_{\omega}$ has at least four-dimensional generalized kernel with
the eigenvectors (\ref{kernel}) and two generalized eigenvectors
(see \cite{P05} for details). The eigenvectors of the linearized
operator $\sigma H_{\omega}$ satisfy the $\sigma$-orthogonality
constraints:
\begin{eqnarray}
\label{constraint1} ({\bf u}_0, {\bf U}) & = & \int_{\mathbb{R}}
\left( \bar{u}_0 U_1 + u_0 U_2 + \bar{v}_0 U_3 +
v_0 U_4 \right) dx = 0, \\
\label{constraint2} ({\bf u}_0', \sigma {\bf U}) & = &
\int_{\mathbb{R}} \left( \bar{u}_0' U_1 - u_0' U_2 + \bar{v}_0' U_3
- v_0' U_4 \right) dx = 0.
\end{eqnarray}
The constraints (\ref{constraint1}) and (\ref{constraint2})
represent first variations of the conserved quantities $Q$ and $P$
in (\ref{quantity-Q}) and (\ref{quantity-P}) at the linearization
(\ref{linearization}).

It follows from the explicit form of $H_{\omega}$ and from Corollary
\ref{cor-property2} that the eigenvalue problem $H_{\omega} {\bf U}
= \mu {\bf U}$ has two reductions:
\begin{equation}
\label{reduction12} (i) \; U_1 = U_4, \; U_2 = U_3, \qquad (ii) \;
U_1 = -U_4, \; U_2 = -U_3.
\end{equation}
Our main result on the block-diagonalization of the energy operator
$H_{\omega}$ and the linearized coupled-mode system
(\ref{linearized-Ham-system}) is based on the reductions
(\ref{reduction12}).

\begin{thm}
\label{theorem-diagonalization} Let $W \in \mathbb{R}$ satisfy
properties (1)-(3). Let $(u_0,v_0)$ be a decaying solution of the
system (\ref{ODE-system}) on $x \in \mathbb{R}$, where $v_0 =
\bar{u}_0$. There exists an orthogonal similarity transformation
$S$, such that $S^{-1} = S^T$, where
$$
S = \frac{1}{\sqrt{2}} \left (
\begin{array}{cccc}1 & 0 & 1 & 0 \\ 0 & 1 & 0 & 1 \\ 0 & 1 & 0 & -1 \\
1 & 0 & -1 & 0 \end{array}\right),
$$
that simultaneously block-diagonalizes the energy operator
$H_{\omega}$,
\begin{equation}
\label{block1}
S^{-1} H_{\omega} S = \left ( \begin{array}{cc} H_+ & 0 \\
0 & H_- \end{array}\right) \equiv H,
\end{equation}
and the linearized operator $\sigma H_{\omega}$
\begin{equation}
\label{block2} S^{-1} \sigma H_{\omega} S = \sigma \left(
\begin{array}{cc} 0 & H_- \\ H_+ & 0 \end{array}\right) \equiv i L,
\end{equation}
where $H_{\pm}$ are two-by-two Dirac operators:
\begin{equation}
\label{dirac1}
H_{\pm} = \left ( \begin{array}{cc} \omega - i \partial_x  & \mp 1 \\
\mp 1 & \omega + i \partial_x \end{array}\right) + V_{\pm}(x),
\end{equation}
and
\begin{equation}
\label{dirac2} V_{\pm} = \left(
\begin{array}{cc} \partial_{\bar{u}_0 u_0}^2 \pm \partial_{\bar{u}_0
\bar{v}_0}^2 & \partial_{\bar{u}_0^2}^2 \pm \partial_{\bar{u}_0 v_0}^2 \\
\partial_{u_0^2}^2 \pm \partial_{u_0 \bar{v}_0}^2 &
\partial_{\bar{u}_0 u_0}^2 \pm \partial_{u_0 v_0}^2
\end{array}\right) W(u_0,\bar{u}_0,v_0,\bar{v}_0).
\end{equation}
\end{thm}

\begin{proof}
Applying the similarity transformation to the operator
$D(\partial_x)$ in (\ref{operator-D}), we have the first terms in
Dirac operators $H_{\pm}$. Applying the same transformation to the
potential $V(x)$ in (\ref{operator-V}) and using Corollary
\ref{cor-property2}, we have the second term in the Dirac
operators $H_{\pm}$. The same transformation is applied similarly
to the linearized operator $\sigma H_{\omega}$ with the result
(\ref{block2}).
\end{proof}

\begin{cor}
The linearized coupled-mode system (\ref{linearized-Ham-system}) is
equivalent to the block-diagonalized eigenvalue problems
\begin{equation}
\label{eigenvalue} \sigma_3 H_- \sigma_3 H_+ {\bf V}_1 = \gamma {\bf
V}_1, \qquad \sigma_3 H_+ \sigma_3 H_- {\bf V}_2 = \gamma {\bf V}_2,
\qquad \gamma = - \lambda^2,
\end{equation}
where ${\bf V}_{1,2} \in \mathbb{C}^2$ and $\sigma_3$ is the Pauli's
diagonal matrix of $(1,-1)$.
\end{cor}

\begin{cor}
Let ${\bf u}_0 = (u_0,\bar{u}_0) \in \mathbb{C}^2$ and $({\bf
f},{\bf g})$ be a standard inner product for ${\bf f},{\bf g} \in
L^2(\mathbb{R},\mathbb{C}^2)$. Dirac operators $H_{\pm}$ have
simple kernels with the eigenvectors
\begin{equation}
\label{kernel-reduced} H_+ {\bf u}_0' = 0, \qquad H_- \sigma_3
{\bf u}_0 = 0,
\end{equation}
while the vectors ${\bf V}_{1,2}$ satisfy the constraints
\begin{equation}
\label{constraints-reduced} ({\bf u}_0,{\bf V}_1) = 0, \qquad
({\bf u}_0',\sigma_3 {\bf V}_2 ) = 0.
\end{equation}
\end{cor}

\begin{rem}
Block-diagonalization described in Theorem
\ref{theorem-diagonalization} has nothing in common with the
explicit diagonalization used in reduction (9.2) of \cite{PS03} for
the particular potential function (\ref{example}) with $a_1 = a_2 =
a_4 = 0$ and $a_3 = 1$. Moreover, the reduction (9.2) of \cite{PS03}
does not work for $\omega \neq 0$, while gap solitons do not exist
in this particular model for $\omega = 0$.
\end{rem}

We illustrate block-diagonalization of the eigenvalue problem
(\ref{eigenvalue}) for the quadric potential function
(\ref{example}). When $a_1 = 1$, $a_2 = \rho$ and $a_3 = a_4 = 0$,
the decaying solution $u_0(x)$ is given by
(\ref{soliton-explicit-1}) and the potential matrices $V_{\pm}(x)$
in the Dirac operators $H_{\pm}$ in (\ref{dirac1})--(\ref{dirac2})
are found in the explicit form:
\begin{equation}
\label{Dirac-example1}
V_+ = (1 + \rho) \left( \begin{array}{cc} 2 |u_0|^2 & u_0^2 \\
\bar{u}_0^2  & 2 |u_0|^2 \end{array}\right), \qquad
V_- = \left( \begin{array}{cc} 2 |u_0|^2 & (1-\rho) u_0^2 \\
(1-\rho) \bar{u}_0^2  & 2 |u_0|^2 \end{array}\right).
\end{equation}
When $a_1 = a_2 =0$, $a_3 = 1$ and $a_4 = s$, the decaying solution
$u_0(x)$ is given by either (\ref{soliton-explicit-2a}) or
(\ref{soliton-explicit-2b}) and the potential matrices $V_{\pm}(x)$
take the form:
\begin{equation}
\label{Dirac-example2a}
V_+ = 3 \left( \begin{array}{cc} u_0^2 + \bar{u}_0^2 & 2 |u_0|^2 \\
2 |u_0|^2  & u_0^2 + \bar{u}_0^2 \end{array}\right) +
s \left( \begin{array}{cc} 2 |u_0|^2 & u_0^2 + 3 \bar{u}_0^2 \\
\bar{u}_0^2 + 3 u_0^2  & 2 |u_0|^2 \end{array}\right),
\end{equation}
\begin{equation}
\label{Dirac-example2b}
V_- = \left( \begin{array}{cc} u_0^2 + \bar{u}_0^2 & -2 |u_0|^2 \\
-2 |u_0|^2  & u_0^2 + \bar{u}_0^2 \end{array}\right) +
s \left( \begin{array}{cc} 0 & -u_0^2-\bar{u}_0^2 \\
-u_0^2-\bar{u}_0^2  & 0 \end{array}\right).
\end{equation}
Numerical computations of eigenvalues of the Dirac operators
$H_{\pm}$ and the linearized operator $L$ in (\ref{block1}) and
(\ref{block2}) are developed for the explicit examples
(\ref{Dirac-example1}) and
(\ref{Dirac-example2a})--(\ref{Dirac-example2b}).

\section{Numerical computations of eigenvalues}

Numerical discretization and truncation of the linearized
coupled-mode system (\ref{linearized-Ham-system}) leads to an
eigenvalue problem for large matrices \cite{S01}. Parallel software
libraries were recently developed for computations of large
eigenvalue problems \cite{G01}. We shall use Scalapack library and
distribute computations of eigenvalues of the system
(\ref{linearized-Ham-system}) for different parameter values between
parallel processors of the SHARCnet cluster Idra using Message
Passing Interface \cite{McMaster}.

We implement a numerical discretization of the linearized
coupled-mode system (\ref{linearized-Ham-system}) using the
Chebyshev interpolation method \cite{S02}. The main advantage of the
Chebyshev grid is that clustering of the grid points occurs near the
end points of the interval and this clustering prevents the
appearance of spurious complex eigenvalues from the discretization
of the continuous spectrum. If the eigenvector is analytic in a
strip near the interpolation interval, the corresponding Chebyshev
spectral derivatives converge geometrically, with an asymptotic
convergence factor determined by the size of the largest ellipse in
the domain of analyticity \cite{S02}.

The continuous spectrum for the linearized coupled-mode system
(\ref{linearized-Ham-system}) can be found from the no-potential
case $V(x) \equiv 0$. It consists of two pairs of symmetric branches
on the imaginary axis $\lambda \in i \mathbb{R}$ for $|{\rm
Im}(\lambda)| > 1 - \omega$ and $|{\rm Im}(\lambda)| > 1 + \omega$
\cite{B98,DG04}. In the potential case $V(x) \neq 0$, the continuous
spectrum does not move, but the discrete spectrum appears. The
discrete spectrum is represented by symmetric pairs or quartets of
isolated non-zero eigenvalues and zero eigenvalue of algebraic
multiplicity four for the generalized kernel of $\sigma H_{\omega}$
\cite{B98,DG04}. We note that symmetries of the Chebyshev grid
preserve symmetries of the linearized coupled-mode system
(\ref{linearized-Ham-system}).

We shall study eigenvalues of the energy operator $H_{\omega}$, in
connection to eigenvalues of the linearized operator $\sigma
H_{\omega}$. It is well known \cite{S01,S02} that Hermitian matrices
have condition number one, while non-Hermitian matrices may have
large condition number. As a result, numerical computations for
eigenvalues and eigenvectors have better accuracy and faster
convergence for self-adjoint operators \cite{S01,S02}. We will use
the block-diagonalizations (\ref{block1}) and (\ref{block2}) and
compute eigenvalues of $H_+$, $H_-$, and $L$. The block-diagonalized
matrix can be stored in a special format which requires twice less
memory than a full matrix and it accelerates computations of
eigenvalues approximately in two times.

Figure 1 displays the pattern of eigenvalues and instability
bifurcations for the symmetric quadric potential (\ref{example})
with $a_1 = 1$ and $a_2 = a_3 = a_4 = 0$. The decaying solution
$u_0(x)$ and the potential matrices $V_{\pm}(x)$ are given by
(\ref{soliton-explicit-1}) and (\ref{Dirac-example1}) with $\rho =
0$. Parameter $\omega$ of the decaying solution $u_0(x)$ is defined
in the interval $-1 < \omega < 1$. Six pictures of Fig. 1 shows the
entire spectrum of $L$, $H_+$ and $H_-$ for different values of
$\omega$. (The continuous movie that shows transformation of
eigenvalues when $\omega$ decreases is available as a multi-media
attachment to this article.)

When $\omega$ is close to $1$ (the gap soliton is close to a
small-amplitude sech-soliton), there exists a single non-zero
eigenvalue for $H_+$ and $H_-$ and a single pair of purely imaginary
eigenvalues of $L$ (see subplot (1) on Fig. 1). The first set of
arrays on the subplot (1) indicates that the pair of eigenvalues of
$L$ becomes visible at the same value of $\omega$ as the eigenvalue
of $H_+$. This correlation between eigenvalues of $L$ and $H_+$ can
be traced throughout the entire parameter domain on the subplots
(1)--(6).

When $\omega$ decreases, the operator $H_-$ acquires another
non-zero eigenvalue by means of the edge bifurcation \cite{KS02},
with no changes in the number of isolated eigenvalues of $L$ (see
subplot (2)). The first complex instability occurs near $\omega
\approx -0.18$, when the pair of purely imaginary eigenvalues of $L$
collides with the continuous spectrum and emerge as a quartet of
complex eigenvalues, with no changes in the number of isolated
eigenvalues for $H_+$ and $H_-$ (see subplot (3)).

The second complex instability occurs at $\omega \approx -0.54$,
when the operator $H_-$ acquires a third non-zero eigenvalue and the
linearized operator $L$ acquires another quartet of complex
eigenvalues (see subplot (4)). The second set of arrays on the
subplots (4)--(6) indicates a correlation between these eigenvalues
of $L$ and $H_-$.

When $\omega$ decreases further, the operators $H_+$ and $H_-$
acquires one more isolated eigenvalue, with no change in the
spectrum of $L$ (see subplot (5)). Finally, when $\omega$ is close
to $-1$ (the gap soliton is close to the large-amplitude algebraic
soliton), the third complex instability occurs, correlated with
another edge bifurcation in the operator $H_-$ (see subplot (6)).
The third set of arrays on subplot (6) indicates this correlation.
The third complex instability was missed in the previous numerical
studies of the same system \cite{B98,DG04}. In a narrow domain near
$\omega = -1$, the operator $H_+$ has two non-zero eigenvalues, the
operator $H_-$ has five non-zero eigenvalues and the operator $L$
has three quartets of complex eigenvalues.

Figure 2 displays the pattern of eigenvalues and instability
bifurcations for the symmetric quadric potential (\ref{example})
with $a_1 = a_2 = a_4 = 0$ and $a_3 = 1$. The decaying solution
$u_0(x)$ and the potential matrices $V_{\pm}(x)$ are given by
(\ref{soliton-explicit-2a}) and (\ref{Dirac-example2a}) with $\omega
> 0$ and $s = 0$. Eigenvalues in the other case $\omega < 0$ can be
found from those in the case $\omega > 0$ by reflections.

When $\omega$ is close to $1$ (the gap soliton is close to a
small-amplitude sech-soliton), there exists one non-zero eigenvalue
of $H_-$ and no non-zero eigenvalues of $L$ and $H_+$ (see subplot
(1)). When $\omega$ decreases, two more non-zero eigenvalues
bifurcate in $H_-$ from the left and right branches of the
continuous spectrum, with no change in non-zero eigenvalues of $L$
(see subplot (2)). The first complex bifurcation occurs at $\omega
\approx 0.45$, when a quartet of complex eigenvalues occurs in $L$,
in correlation with two symmetric edge bifurcations of $H_+$ from
the left and right branches of the continuous spectrum (see subplot
(3)). The first and only set of arrays on the subplots (3)-(6)
indicates a correlation between eigenvalues of $L$ and $H_+$, which
is traced through the remaining parameter domain of $\omega$. The
inverse complex bifurcation occurs at $\omega \approx 0.15$, when
the quartet of complex eigenvalues merge at the edge of the
continuous spectrum into a pair of purely imaginary eigenvalues (see
subplot (5)). No new eigenvalue emerge for smaller values of
$\omega$. When $\omega$ is close to $0$ (the gap soliton is close to
the non-decaying solution), the operator $H_+$ has two non-zero
eigenvalues, the operator $H_-$ has three non-zero eigenvalues and
the operator $L$ has one pair of purely imaginary eigenvalues (see
subplot (6)).

We mention two other limiting cases of the symmetric quadric
potential (\ref{example}). When $a_1 = a_3 = a_4 = 0$ and $a_2 = 1$,
the coupled-mode system is an integrable model and no non-zero
eigenvalues of $L$ exist, according to the exact solution of the
linearization problem \cite{KL96,KL97}. When $a_1 = a_2 = a_3 = 0$
and $a_4 = \pm 1$, one branch of decaying solutions $u_0(x)$ exists
for either sign, according to (\ref{soliton-explicit-2a}) and
(\ref{soliton-explicit-2b}). The pattern of eigenvalues and
instability bifurcations repeats that of Fig. 2.

Numerical results reported above imply that the number of isolated
non-zero eigenvalues of the linearized operator $L$ is bounded from
above by the total number of non-zero isolated eigenvalues of the
energy operators $H_+$ and $H_-$. Furthermore, there exists a
correlation between edge bifurcations in the operator $L$ and those
in the Dirac operators $H_+$ and $H_-$. These analytical questions
will be addressed in the future work.

\appendix
\section{Conditions for existence of gap solitons in the homogeneous potential function}

We shall consider the homogeneous potential function $W \in
\mathbb{R}$ of the monomial order $2n$ that satisfies properties
(1)-(3). The general representation of $W(u,\bar{u},v,\bar{v})$ is
given by
\begin{equation}
\label{hexample} W = \sum_{s=0}^{n} \sum_{k=0}^{n-s} a_{k,s} \left(
u^s \bar{v}^s + \bar{u}^s v^s \right) |u|^{2n-2k-2s} |v|^{2k},
\end{equation}
where $a_{k,s}$ are real-valued coefficients which are subject to
the symmetry conditions: $a_{k_1,s} = a_{k_2,s}$ if $k_1 + k_2= n-s$
for $s = 0,1,...,n-1$. Let's introduce new parameters
$$
A_s = \sum_{k=0}^{n-s} a_{k,s}, \qquad s = 0,1,...,n.
$$
Using the variables $(Q,\Theta)$ defined in
(\ref{representation-soliton}) with $\Phi(x) = \Phi_0 = 0$, we
rewrite the ODE system (\ref{ODE-Hamiltonian}) in the explicit form:
\begin{equation}
\label{hexplicitform} \left \{
\begin{array}{c}
Q' = 2 Q \sin(2 \Theta) - 2 Q^n \sum_{s=0}^{n} s A_s \sin(2 s \Theta) \\
\Theta' = - \omega + \cos(2 \Theta) - n Q^{n-1} \sum_{s=0}^{n} A_s
\cos(2 s \Theta)
\end{array}
\right.
\end{equation}
There exists a first integral of the system (\ref{hexplicitform}):
$$
-\omega Q + \cos(2 \Theta) Q - Q^n \sum_{s=0} ^{n} A_s \cos(2 s
\Theta) = C_0,
$$
where $C_0 = 0$ from the zero boundary conditions $Q(x) \to 0$ as
$|x| \to \infty$. As a result, the second-order system
(\ref{hexplicitform}) is reduced to the first-order ODE
\begin{eqnarray}
\label{hcoupled-mode-system-appendix} \Theta'(x) = (n-1)(\omega -
\cos(2 \Theta)),
\end{eqnarray}
while the function $Q(x) \geq 0$ can be found from $\Theta(x)$ as
follows:
\begin{equation}
\label{hQ} Q^{n-1} = \frac{(\cos(2\Theta) - \omega)}{\sum_{s=0}^{n}
A_s \cos(2 s \Theta)}.
\end{equation}

We consider the quadric potential function $W$ given by
(\ref{example}). Using (\ref{hcoupled-mode-system-appendix}) for the
case $n=2$ we obtain:
\begin{eqnarray}
\label{coupled-mode-system-appendix} \Theta'(x) = \omega - \cos(2
\Theta),
\end{eqnarray}
and the correspondence:
$$
A_0 = \frac{a_1 + a_2 + a_4}{2}, \quad A_1 = 2 a_3, \quad A_2 =
\frac{a_4}{2}.
$$
We rewrite the representation (\ref{hQ}) for $Q(x)$ as follows:
\begin{equation}
Q = \frac{(t - \omega)}{\phi(t)}; \qquad Q \geq 0
\end{equation}
where
$$
t = \cos(2 \Theta), \qquad \phi(t) = a_4 t^2 + 2 a_3 t + \frac{a_1 +
a_2}{2},
$$
such that $t \in [-1,1]$. Let's consider two cases:
\begin{equation}
\label{twocases} \left \{
\begin{array}{c}
t \geq \omega; \quad \phi(t) \geq 0 \quad \Rightarrow Q^+\\
t \leq \omega; \quad \phi(t) \leq 0 \quad \Rightarrow Q^-\\
\end{array}
\right.
\end{equation}
We can solve the first-order ODE
(\ref{coupled-mode-system-appendix}) using the substitution $z =
\tan(\Theta)$, such that
$$
t = \frac{1 - z^2}{1 + z^2} \qquad z^2 = \frac{1-t}{1+t}.
$$
After integration with the symmetry constraint $\Theta(0) = 0$, we
obtain the solution
\begin{equation}
\left|\frac{(z - \sqrt{\mu})}{( z + \sqrt{\mu})}\right| = e^{2 \beta
x},
\end{equation}
where
$$
\beta = \sqrt{1 - \omega^2}, \qquad \mu = \frac{1-\omega}{1+\omega}
$$
and $-1 < \omega < 1$. Two separate cases are considered:
\begin{equation}
\label{Qcase1} |z| \leq \sqrt{\mu} \qquad z = - \sqrt{\mu}
\frac{\sinh(\beta x)}{\cosh(\beta x)} \qquad t = \frac{\cosh^2(\beta
x) - \mu \sinh^2(\beta x)}{\cosh^2(\beta x) + \mu \sinh^2(\beta x)},
\end{equation}
where $t \geq \omega$, and
\begin{equation}
\label{Qcase2}
 |z| \geq \sqrt{\mu} \qquad z = - \sqrt{\mu}
\frac{\cosh(\beta x)}{\sinh(\beta x)} \qquad t = \frac{\sinh^2(\beta
x) - \mu \cosh^2(\beta x)}{\sinh^2(\beta x) + \mu \cosh^2(\beta x)},
\end{equation}
where $t \leq \omega$. Let's introduce new parameters
\begin{eqnarray*}
A & = & -2a_3 + a_4 + \frac{a_1+a_2}{2}, \\
B & = & -2a_4 + a_1 + a_2, \\
C & = & 2a_3 + a_4 + \frac{a_1+a_2}{2}.
\end{eqnarray*}
It is clear that $A = \phi(-1)$ and $C = \phi(1)$. If $t \geq
\omega$ and $\phi(t) \geq 0$, it follows from (\ref{twocases}) and
(\ref{Qcase1}) that
\begin{equation}
Q^+(x) = \frac {(1-\omega)((\mu+1) \cosh^2(\beta x) - \mu)}{(A \mu^2
+ B \mu + C) \cosh^4( \beta x) - (B\mu + 2A \mu^2 )\cosh^2 (\beta x)
 + A \mu^2}.
\end{equation}
If $t \leq \omega $ and  $ \phi(t) \leq 0$, it follows from
(\ref{twocases}) and (\ref{Qcase2}) that
\begin{equation}
Q^-(x) = \frac {(\omega - 1)((\mu+1) \cosh^2(\beta x) - 1)}{(A \mu^2
+ B \mu + C) \cosh^4( \beta x) - (B\mu + 2C )\cosh^2 (\beta x)
 + C}.
\end{equation}
The asymptotic behavior of the function $Q(x)$ at infinity depends
on the location of the zeros of the function $\psi(\mu) = A \mu^2 +
B \mu + C $. The function $\psi(\mu)$ is related to the function
$\phi(t)$, e.g. if $\psi(\mu) = 0$ then $\phi(\omega) = 0$.

\subsection{Case $A < 0$, $C > 0$}
In this case the quadratic polynomial $\phi(t)$ has exactly one root
$\phi(t_1)=0$ such that $t_1 \in (-1,1)$. We have two branches of
decaying solutions with the positive amplitude $Q(x)$. One branch
occurs for $t_1 < \omega \leq 1$ with $Q(x) = Q^+(x)$ and the other
one occurs for $ -1 \leq \omega < t_1$ with $Q(x) = Q^-(x)$. At the
point $\omega = t_1$, the solution is bounded and decaying.

\subsection{Case $A > 0$, $C > 0$}
In this case the quadratic polynomial $\phi(t)$ has no roots or has
exactly two roots on $(-1,1)$. If $\phi(t)$ does not have any roots
on $(-1,1)$, we have a decaying solution with the positive amplitude
$Q(x)$ for any $-1 < \omega < 1$ with $Q(x) = Q^+(x)$. If $\phi(t)$
has two roots $\phi(t_1) =0$ and $\phi(t_2)=0$ such that $t_1,t_2
\in (-1,1)$ then we have a decaying solution with $Q(x) = Q^+(x)$
only on the interval $\max(t_1, t_2)< \omega \leq 1 $. At the point
$\omega = \max(t_1, t_2)$, the solution becomes bounded but
non-decaying if $t_1 \neq t_2$ and unbounded if $t_1 = t_2$.

\subsection{Case $A < 0$, $C < 0$}
In this case the quadratic polynomial $\phi(t)$ has no roots or has
exactly two roots on $(-1,1)$. If $\phi(t)$ does not have any roots
on $(-1,1)$, we have a decaying solution with the positive amplitude
$Q(x)$ for any $-1 < \omega < 1$ with $Q(x) = Q^-(x)$. If $\phi(t)$
has two roots $\phi(t_1) =0$ and $\phi(t_2)=0$ such that $t_1,t_2
\in (-1,1)$ then we have a decaying solution with $Q(x) = Q^-(x)$
only on the interval $-1 \leq \omega < \min(t_1, t_2)$. At the point
$\omega = \min(t_1, t_2)$, the solution becomes bounded but
non-decaying if $t_1 \neq t_2$ and unbounded if $t_1 = t_2$.

\subsection{Case $A > 0$, $C < 0$}
In this case no decaying solutions with positive amplitude $Q(x)$
exist.

\subsection{Special cases}
Two special cases occur when $\phi(1) = 0$ or $\phi(-1) = 0$. If
$\phi(1) = 0$ then $Q^+(x)$ has a singularity at $x = 0$ for any $-1
< \omega < 1$. If $\phi(-1) = 0$ then $Q^-(x)$ has a singularity at
$x = 0$ for any $-1 < \omega < 1$.

\clearpage

\begin{figure}
\centering
\begin{tabular}{cc}
\\

\begin{minipage}{2.5 in}
\includegraphics [width = 6cm] {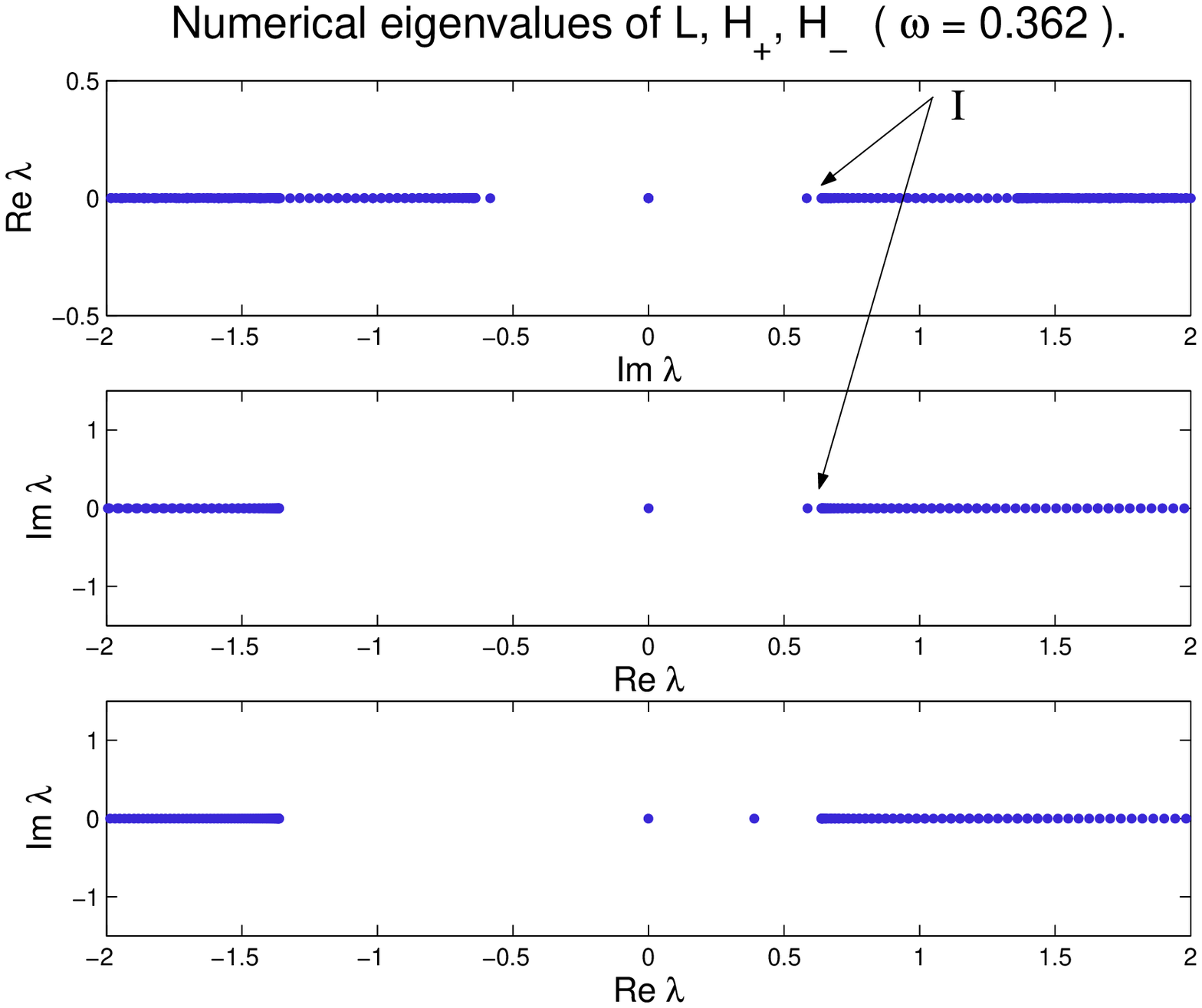}
\end{minipage}
& \begin{minipage}{2.5 in}
\includegraphics [width = 6cm]{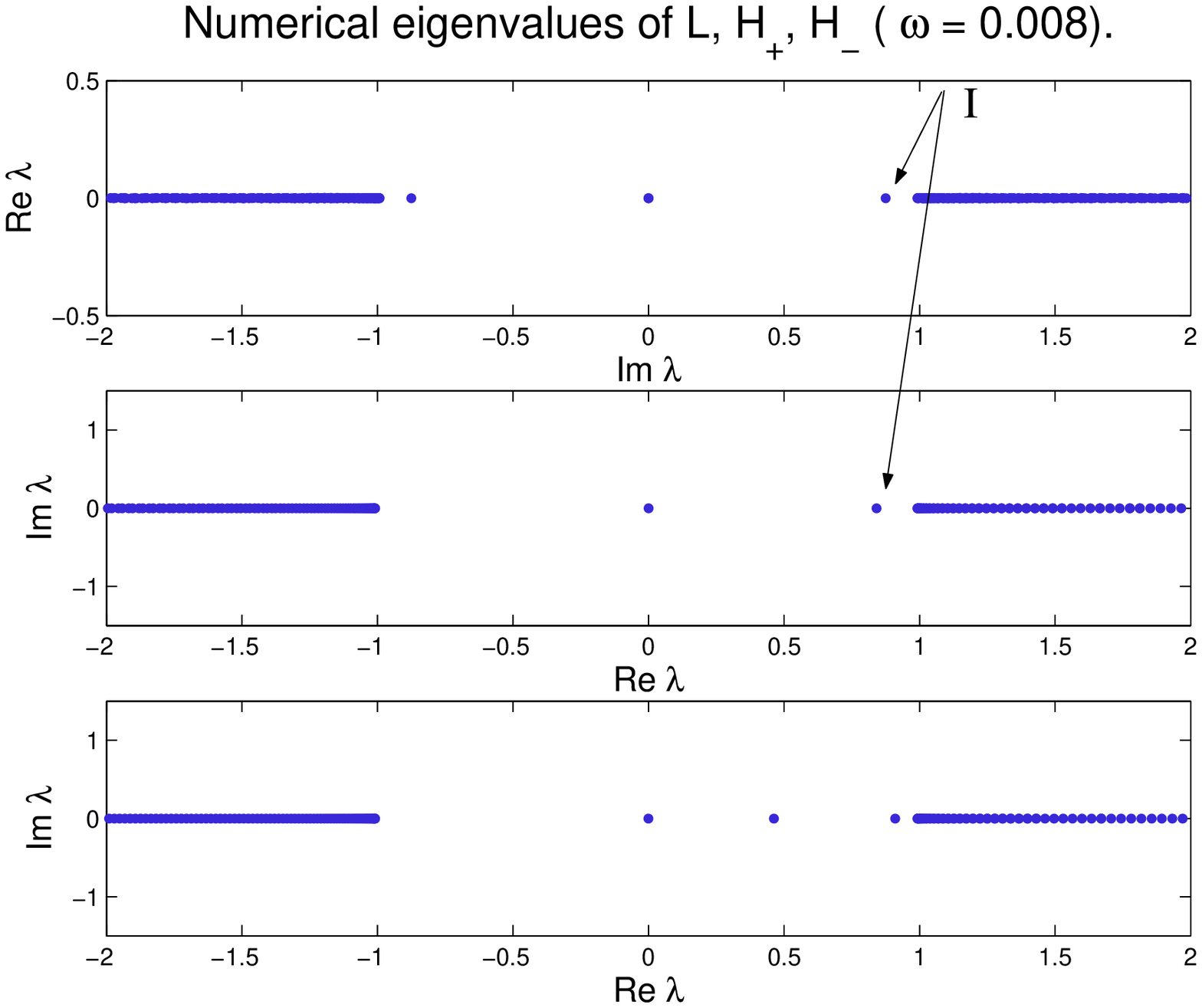}
\end{minipage} \\

\\
\begin{minipage}{2.5 in}
\includegraphics [width = 6cm] {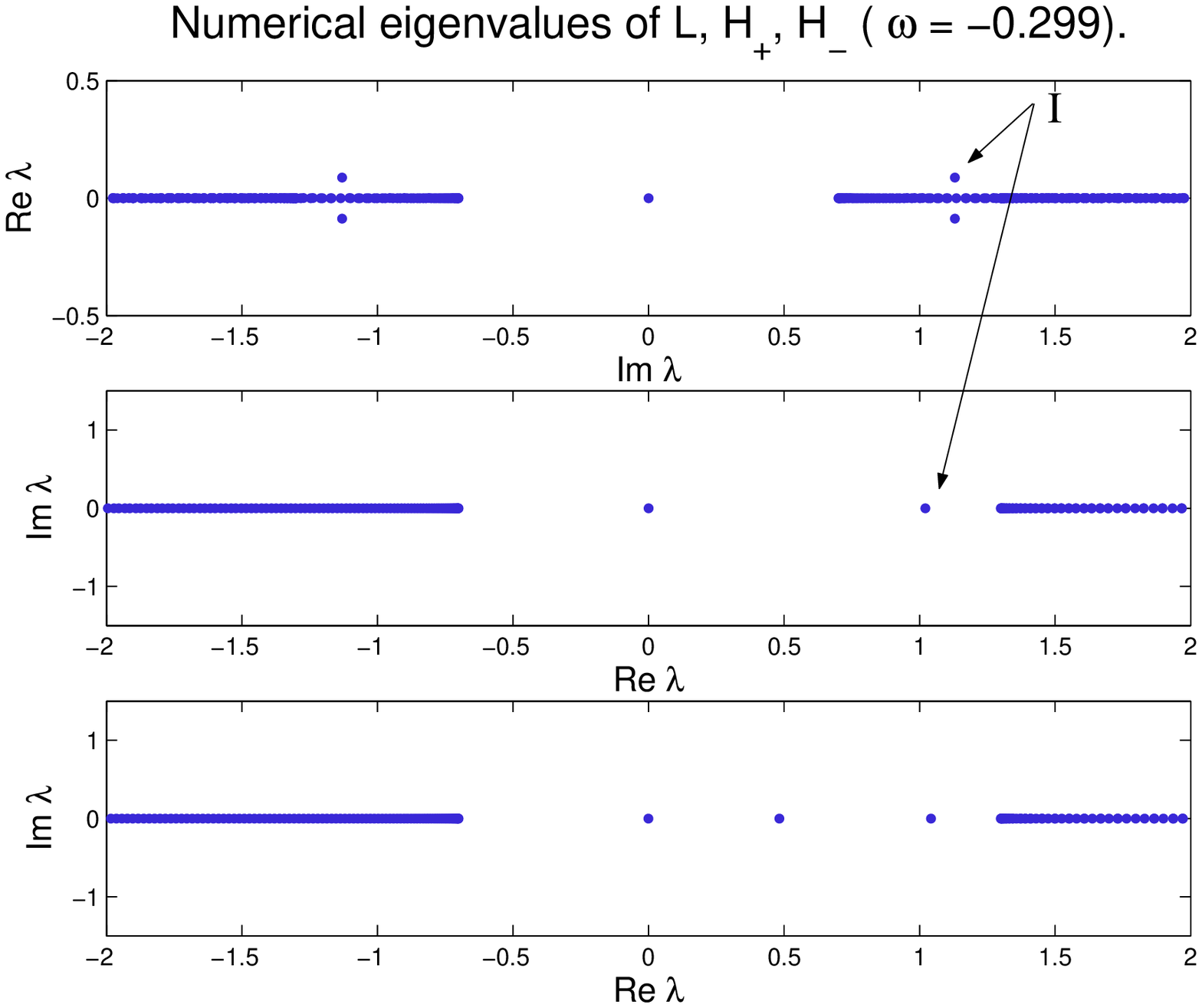}
\end{minipage}
& \begin{minipage}{2.5 in}
\includegraphics [width = 6cm]{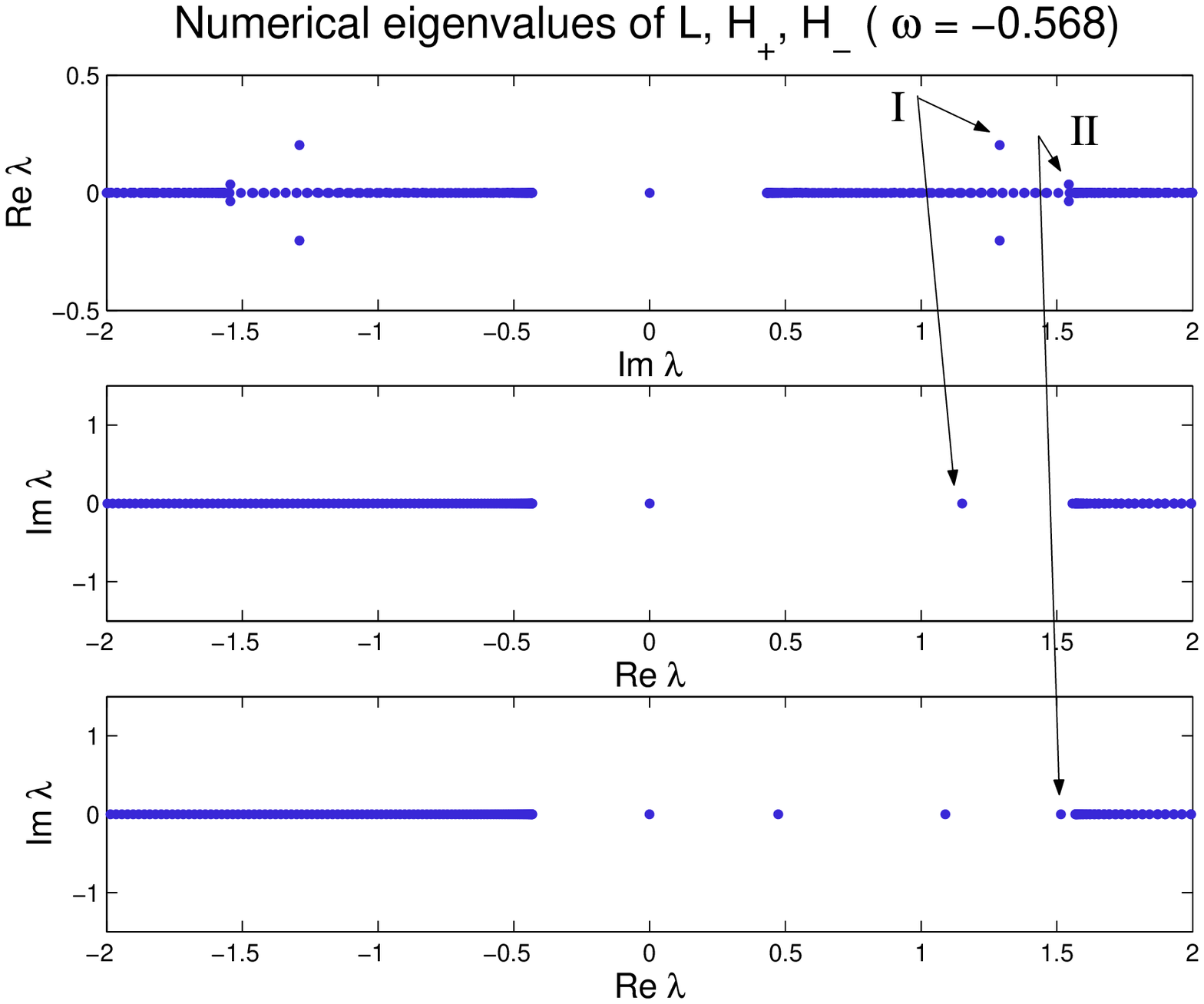}
\end{minipage} \\

\\
\begin{minipage}{2.5 in}
\includegraphics [width = 6cm] {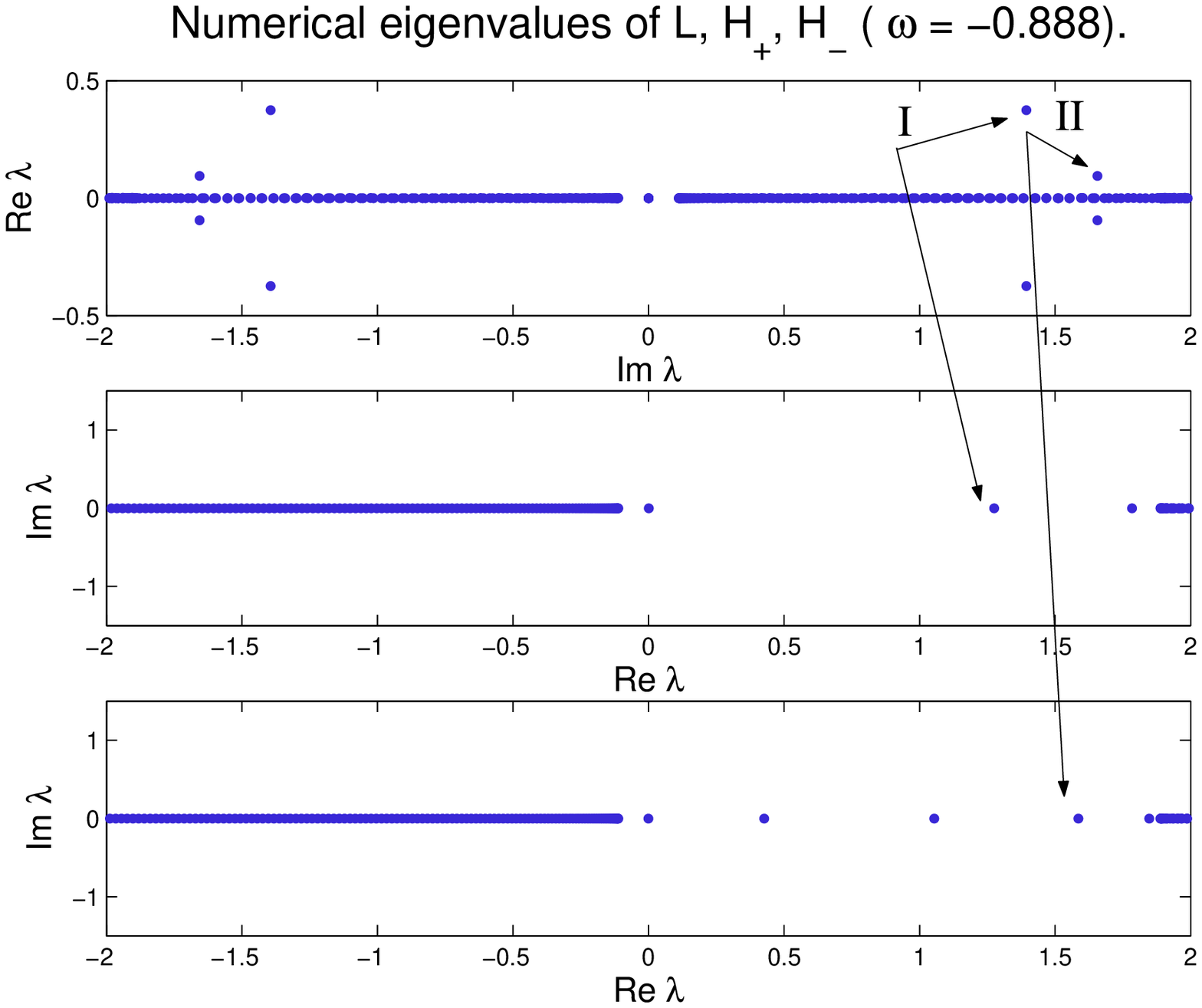}
\end{minipage}
& \begin{minipage}{2.5 in}
\includegraphics [width = 6cm]{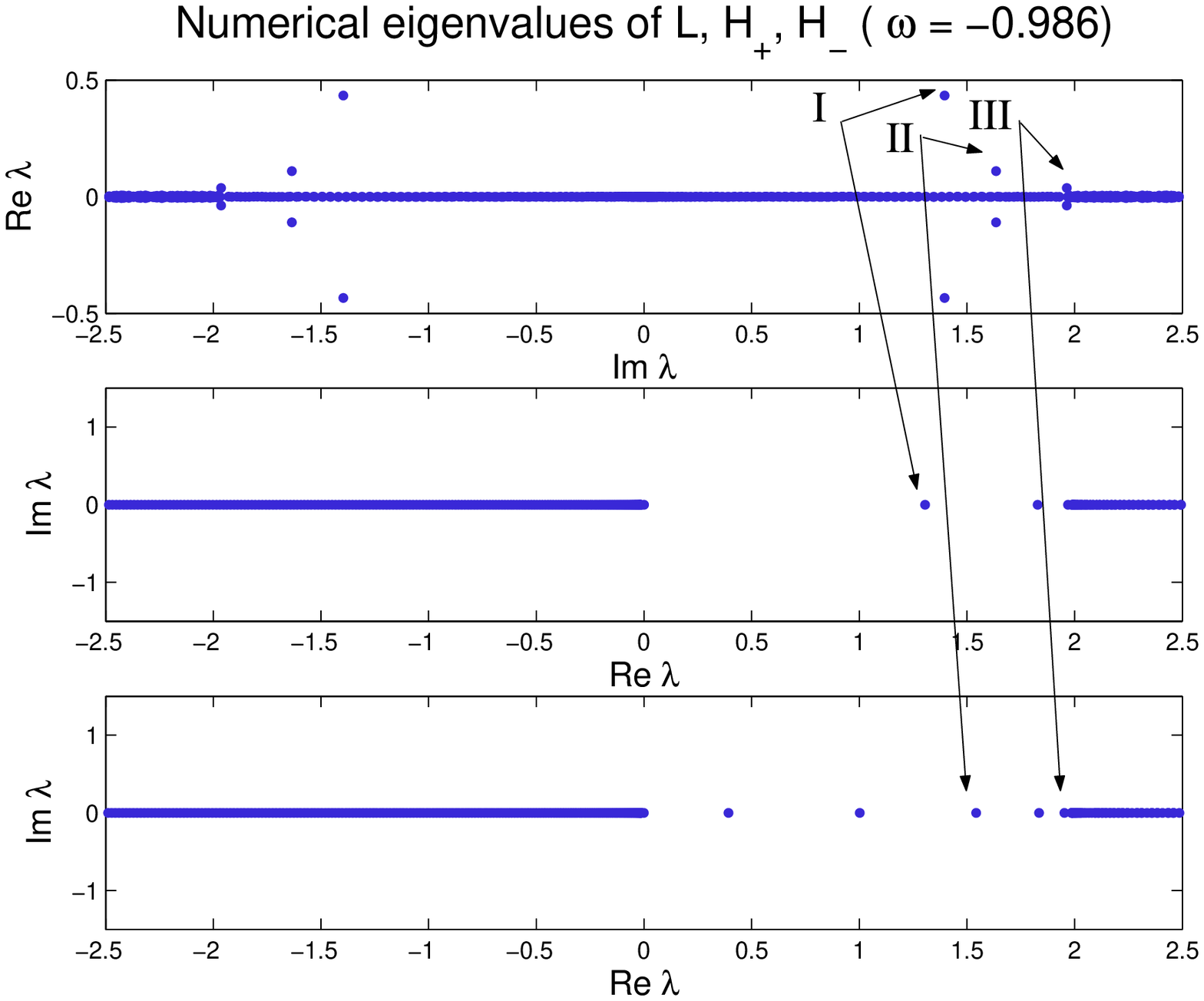}
\end{minipage} \\
\end{tabular}
\end{figure}

\begin{figure}
\caption{Eigenvalues and instability bifurcations for the symmetric
quadric potential (\ref{example}) with $a_1 = 1$ and $a_2 = a_3 =
a_4 = 0$.}
\end{figure}
\clearpage

\begin{figure}
\centering
\begin{tabular}{cc}
\\

\begin{minipage}{2.5 in}
\includegraphics [width = 6cm] {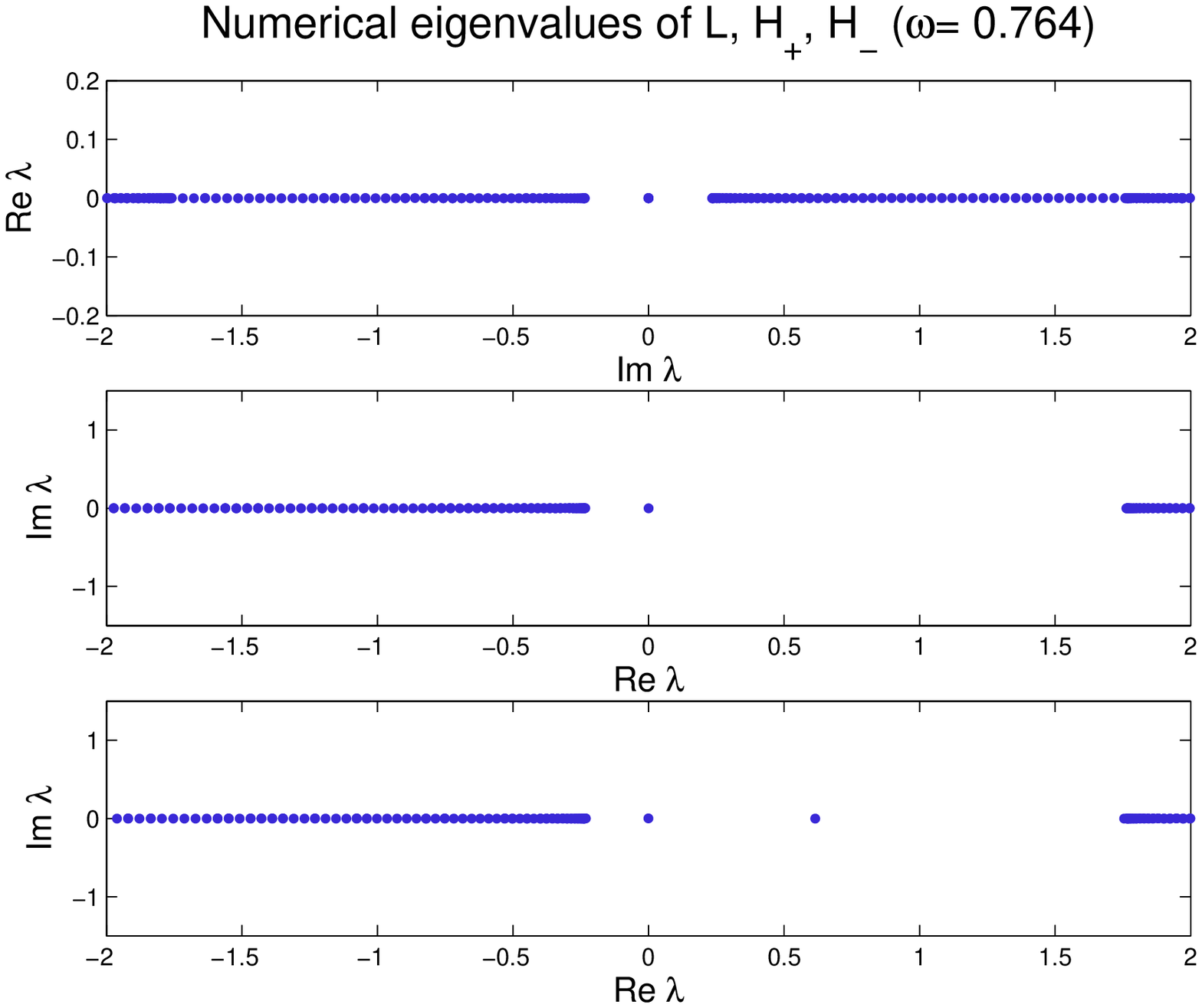}
\end{minipage}
& \begin{minipage}{2.5 in}
\includegraphics [width = 6cm]{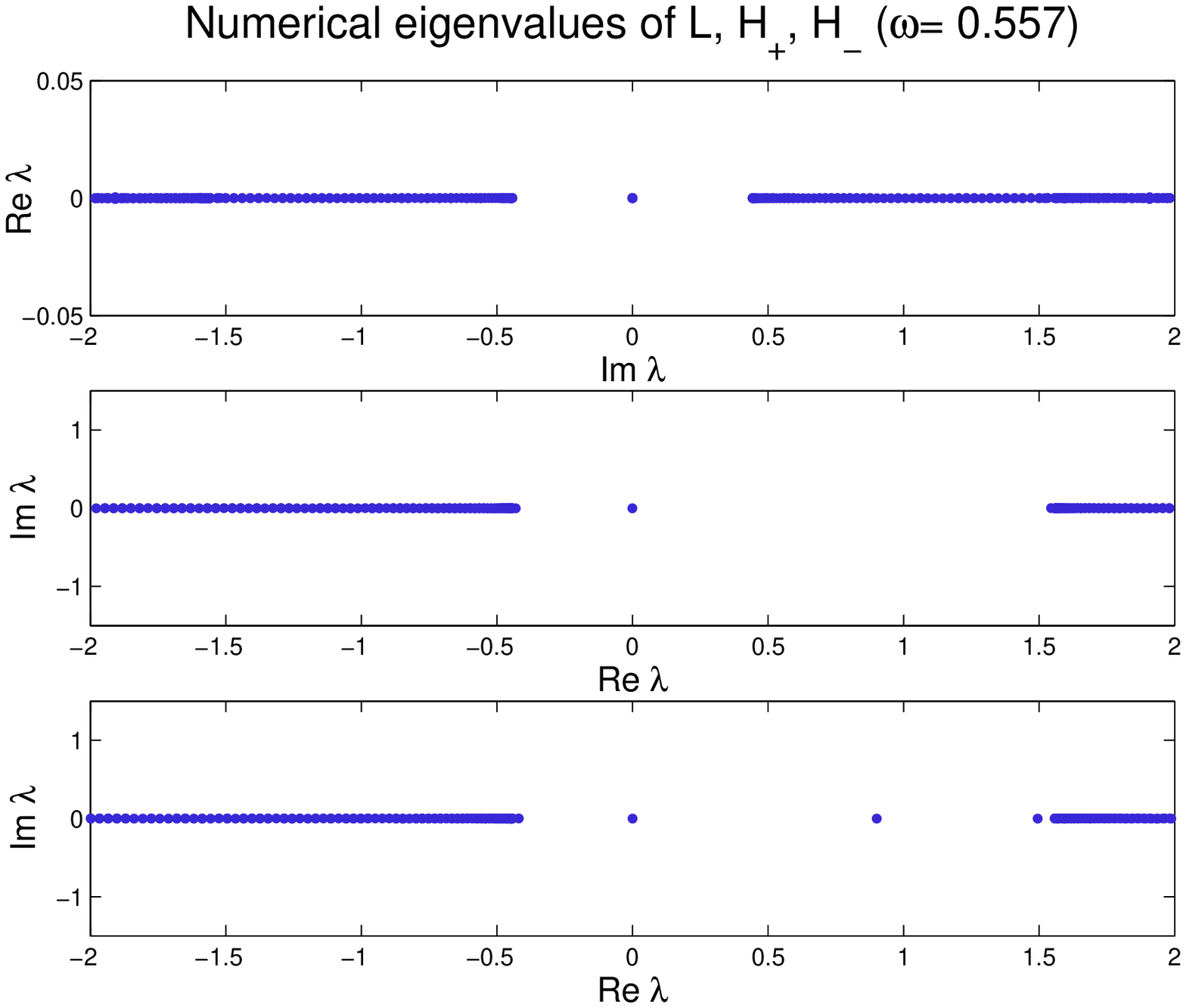}
\end{minipage} \\

\\
\begin{minipage}{2.5 in}
\includegraphics [width = 6cm] {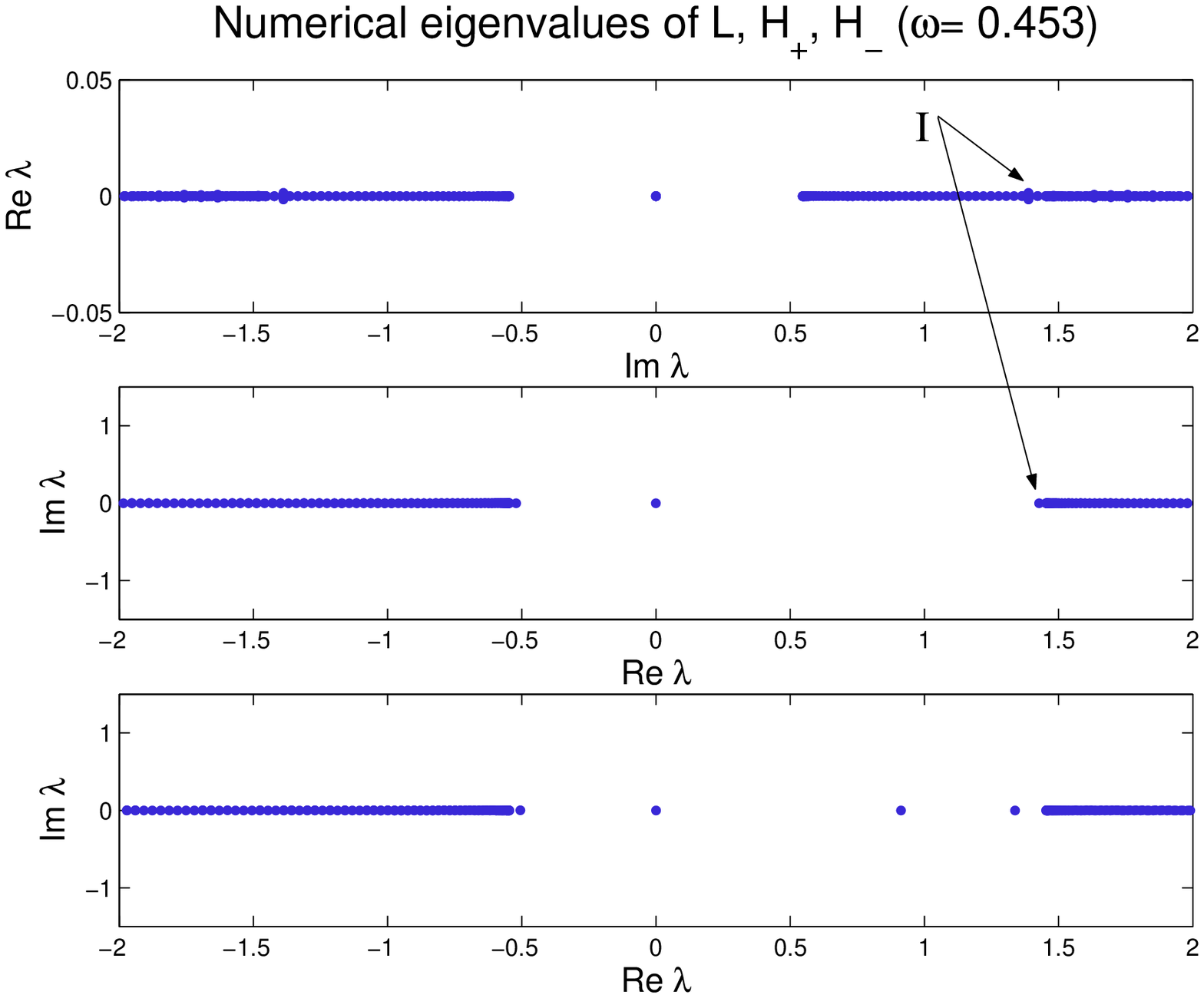}
\end{minipage}
& \begin{minipage}{2.5 in}
\includegraphics [width = 6cm]{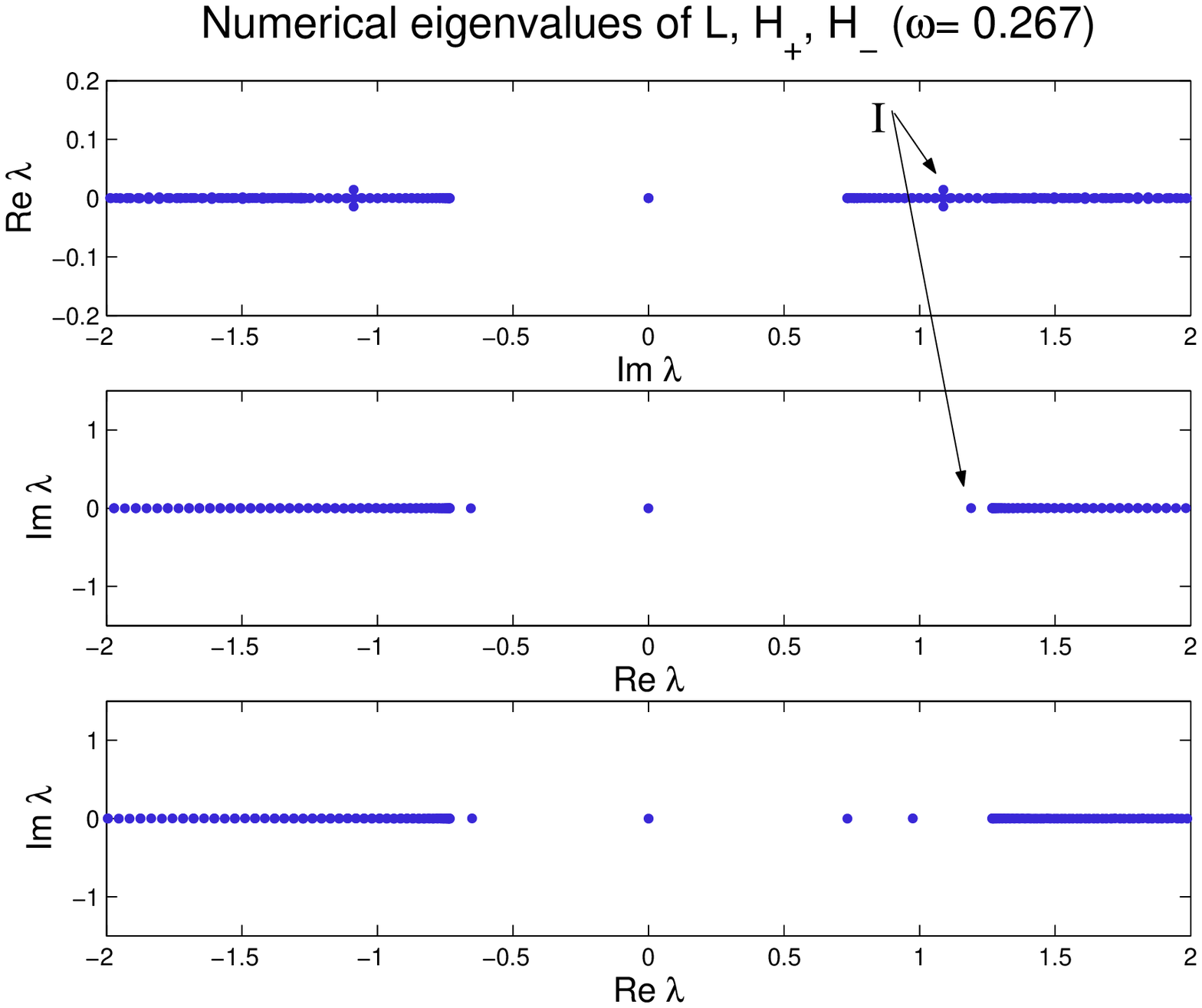}
\end{minipage} \\

\\
\begin{minipage}{2.5 in}
\includegraphics [width = 6cm] {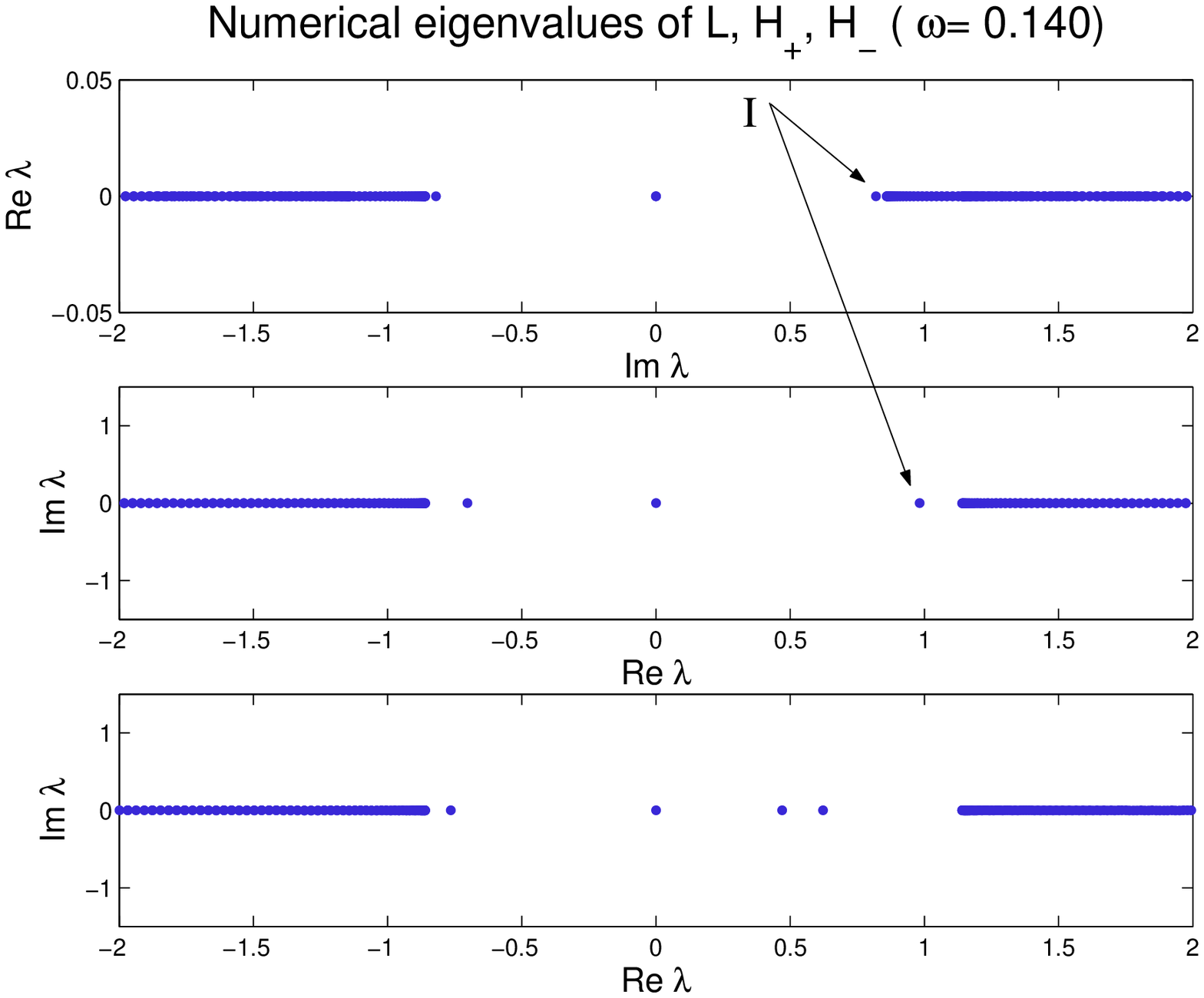}
\end{minipage}
& \begin{minipage}{2.5 in}
\includegraphics [width = 6cm]{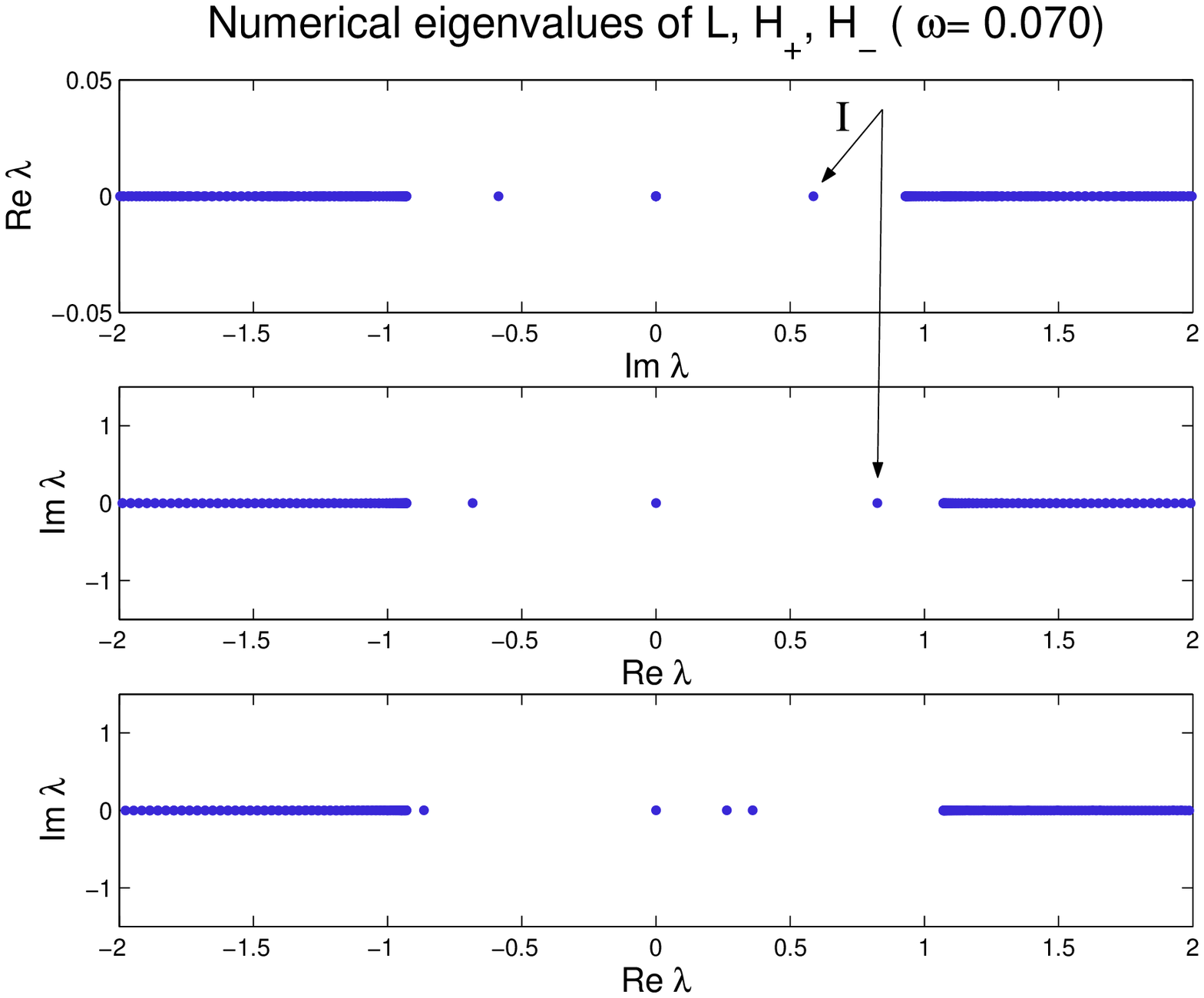}
\end{minipage} \\
\end{tabular}
\end{figure}

\begin{figure}
\caption{Eigenvalues and instability bifurcations for the symmetric
quadric potential (\ref{example}) with $a_3 = 1$ and $a_1 = a_2 =
a_4 = 0$.}
\end{figure}
\clearpage

\bibliographystyle{amsplain}

\begin{thebibliography}{200}
\bibitem{McMaster} Cluser Idra is a part of the SHARCnet network of
parallel processors distributed between eight universities in
Southern Ontario, including McMaster University.
\end{thebibliography}

\end{document}